\documentclass{amsart}
\usepackage[all]{xy}

\newtheorem{theorem}{Theorem}[section]
\newtheorem{lemma}{Lemma}[section]
\newtheorem{proposition}{Proposition}[section]

\newtheorem{definition}{Definition}

\newtheorem{nonexample}{Counterexample}

\theoremstyle{remark}
\newtheorem{remark}{Remark}[section]
\numberwithin{equation}{section}
\numberwithin{lemma}{section}
\numberwithin{proposition}{section}
\numberwithin{definition}{section}

\newcommand{\internalcomment}[1]{}

\begin{document}

\title[A new Quillen model]{A new Quillen model
  for the Morita homotopy theory of DG categories}
\author{Gon{\c c}alo Tabuada}
\address{Universit{\'e} Paris 7 - Denis Diderot, UMR 7586
  du CNRS, case 7012, 2 Place Jussieu, 75251 Paris cedex 05, France.}

\thanks{Supported by FCT-Portugal, scholarship {\tt SFRH/BD/14035/2003}.}

\email{
\begin{minipage}[t]{5cm}
tabuada@math.jussieu.fr
\end{minipage}
}

\begin{abstract}
We construct a new Quillen model, based on the notions of Drinfeld's
DG quotient, \cite{Drinfeld}, and localization pair, for the Morita homotopy theory of
DG categories. This new Quillen model carries a natural closed
symmetric monoidal structure and allow us to re-interpret To{\"e}n's
construction of the internal $\mathsf{Hom}$-functor for the homotopy
category of DG categories as a total right derived internal $\mathsf{Hom}$-functor.
\end{abstract}

\maketitle

\tableofcontents

\section{Introduction}
In this article we propose a solution to the following problem stated by
To{\"e}n in \cite{Toen}, where we suppose that the commutative ground ring
$k$ is a field:
\bigskip

{\it The model category $\mathsf{dgcat}$ together with the symmetric
  monoidal structure $-\otimes-$ is not a symmetric monoidal model
  category, as the tensor product of two cofibrant objects in
  $\mathsf{dgcat}$ is not cofibrant in general. A direct consequence
  of this fact is that the internal $\mathsf{Hom}$ object between
  cofibrant-fibrant objects in $\mathsf{dgcat}$ can not be invariant
  by quasi-equivalences, and thus does not provide internal
  $\mathsf{Hom}$'s for the homotopy categories
  $\mathsf{Ho}(\mathsf{dgcat})$.}
\bigskip

In \cite{Toen}, To{\"e}n has constructed the internal $\mathsf{Hom}$-functor
$\mathsf{rep}_{dg}(-,-)$ for $\mathsf{Ho}(\mathsf{dgcat})$, using a
certain dg category of right quasi-representable bimodules.

In \cite{Drinfeld}, Drinfeld has given an explicit construction of the
dg quotient of a dg category $\mathcal{A}$ modulo a full dg subcategory
$\mathcal{B}$, under certain flatness assumptions that are satisfied
if one works over a field.

In this article, we construct a Quillen model structure on the
category $\mathsf{Lp}$ of localization pairs using Drinfeld's explicit
dg quotient construction. We show that this new model is Quillen
equivalent to the one constructed in \cite{IMRN}\cite{erratum}, and
carries a natural closed symmetric monoidal structure. The tensor
product and internal $\mathsf{Hom}$-functor in $\mathsf{Lp}$ are shown
to be derivable functors, which correspond, under the equivalence
between $\mathsf{Ho}(\mathsf{Lp})$ and $\mathsf{Ho}(\mathsf{dgcat})$,
to the derived tensor product $-\overset{\mathbb{L}}{\otimes}-$ and
$\mathsf{rep}_{dg}(-,-)$ constructed in \cite{Toen}. In particular we
re-interpret the functor $\mathsf{rep}_{dg}(-,-)$ as a total right derived
internal $\mathsf{Hom}$-functor $\mathcal{R}\mathsf{Hom}(-,-)$ in our
new Quillen model.

\section{Acknowledgments}

This article is part of my Ph.~D. thesis under the supervision of
Prof.~B.~Keller. I deeply thank him for several useful discussions and
generous patience. I am very grateful to B.~To{\"e}n for pointing out an
error in a previous version of this article.

\section{Preliminaries}
In what follows, $k$ will denote a field. The
tensor product $\otimes$ will denote the tensor product over $k$. Let
$\mathsf{Ch}(k)$ denote the category of complexes over $k$. By a
{\it dg category}, we mean a differential graded $k$ category, see
\cite{Drinfeld} \cite{dg-cat} \cite{ICM}. For a dg category
$\mathcal{A}$, we denote by $\mathcal{C}_{dg}(\mathcal{A})$ the dg
category of right $\mathcal{A}$ dg modules and by $\hat{ } : \mathcal{A} \rightarrow
\mathcal{C}_{dg}(\mathcal{A})$ the Yoneda dg functor. We write $\mathsf{dgcat}$
for the
category of small dg categories.
It is proven in \cite{IMRN} \cite{erratum} \cite{cras}, that the
category $\mathsf{dgcat}$ admits a structure of cofibrantly generated
model category whose weak equivalences are the Morita
equivalences defined in \cite{IMRN}\cite{erratum}. Recall that we dispose of an explicit set $I=\{Q, S(n)\}$ of
generating cofibrations and  an explicit set $J = \{R(n), F(n), I_n(k_0, \ldots, k_n),
  L_n(k_0, \ldots, k_n), C\}$ of generating trivial cofibrations.

\section{Homotopy of DG functors}

Let $\mathcal{B}$ be a dg category.

\begin{definition}
Let $P(\mathcal{B})$ be the dg category, see~\cite{Drinfeld}, whose
objects are the closed morphismes of degree zero in $\mathcal{B}$
$$ X \stackrel{f}{\longrightarrow} Y\,,$$
that become invertible in $\mathsf{H}^0(\mathcal{B})$.
We define the complex of morphismes
$$ \mathsf{Hom}_{P(\mathcal{B})}(X\stackrel{f}{\rightarrow}Y,
X'\stackrel{f'}{\rightarrow}Y')$$
as the homotopy pull-back in $\mathsf{Ch}(k)$ of the diagram
$$
\xymatrix{
& \mathsf{Hom}_{\mathcal{B}}(Y,Y') \ar[d]^{f^*} \\
\mathsf{Hom}_{\mathcal{B}}(X,X') \ar[r]^{f'_*} & \mathsf{Hom}_{\mathcal{B}}(X,Y')\,.
}
$$
\end{definition}

\begin{remark}
We dispose of the natural commutative diagram in $\mathsf{dgcat}$
$$
\xymatrix{
\mathcal{B} \ar[rr]^{\Delta} \ar[dr]_i & & \mathcal{B}\times
\mathcal{B} \\
& P(\mathcal{B}) \ar[ur]_{p_0\times p_1} & \,,
}
$$
where $i$ is the dg functor that associates to an object $B$ of
$\mathcal{B}$ the morphism $B \stackrel{Id}{\rightarrow}B$ and $p_0$,
resp. $p_1$, is the dg functor that sends a closed morphism $X
\stackrel{f}{\rightarrow} Y$ to $X$, resp. $Y$.
\end{remark}

\begin{lemma}
The dg category $P(\mathcal{B})$ is a path object for $\mathcal{B}$,
see~\cite{Hirschhornn}, in the Quillen model
structure described in \cite{cras}.
\end{lemma}

\begin{proof}
We prove that the dg functor $i$ is a quasi-equivalence. Clearly
the dg functor $i$ induces a quasi-isomorphism in $\mathsf{Ch}(k)$
$$ \mathsf{Hom}_{\mathcal{B}}(X,Y) \stackrel{\sim}{\longrightarrow}
\mathsf{Hom}_{P(\mathcal{B})}(i(X), i(Y))\,,$$
for every object $X,Y \in \mathcal{B}$.
Remark that the functor $\mathsf{H}^0(i)$ is also essentially surjective.
In fact, let $X \stackrel{f}{\rightarrow}Y$ be an object of
$P(\mathcal{B})$. Consider the following morphism in $P(\mathcal{B})$
from $i(X)$ to $X \stackrel{f}{\rightarrow}Y$,
$$
\xymatrix{
X \ar[r]^{Id} \ar@{=}[d]  \ar@{}[dr]|{h=0} & X \ar[d]^f \\
X \ar[r]^f & Y\,,
}
$$
where $h$ denotes de zero homotopy. Remark that it becomes an
isomorphism in $\mathsf{H}^0(P(\mathcal{B}))$ simply because $f$
becomes an isomorphism in $\mathsf{H}^0(\mathcal{B})$. This proves
that the dg functor $i$ is a quasi-equivalence.
We will now show that the dg functor $p_0 \times p_1$ is a fibration
in the Quillen model structure described in \cite{cras}.
Remark first, that by definition of $P(\mathcal{B})$ the dg functor
$p_0 \times p_1$ induces a surjective morphism in $\mathsf{Ch}(k)$
$$ 
\xymatrix{
\mathsf{Hom}_{P(\mathcal{B})}(X \stackrel{f}{\rightarrow}Y, X'
\stackrel{f'}{\rightarrow}Y') \ar@{->>}[rr]^{p_0\times
  p_1} & & \mathsf{Hom}_{\mathcal{B}}(X,X')\times
\mathsf{Hom}_{\mathcal{B}}(Y,Y')
}
$$
for every object $X \stackrel{f}{\rightarrow}Y$ and
$X'\stackrel{f'}{\rightarrow} Y$ in $P(\mathcal{B})$.
We will now show that contractions lift along the dg functor $P(\mathcal{B})
\stackrel{p_0\times p_1}{\longrightarrow} \mathcal{B}\times
\mathcal{B}$. Let $X \stackrel{f}{\rightarrow}Y$
be an object of $P(\mathcal{B})$. Remark that a contraction of $X
\stackrel{f}{\rightarrow}Y$ in $P(\mathcal{B})$ corresponds exactly to
the following morphisms in $\mathcal{B}$, $c_X \in
\mathsf{Hom}^{-1}_{\mathcal{B}}(X,X)$, $c_Y \in
\mathsf{Hom}^{-1}_{\mathcal{B}}(Y,Y)$ and $h \in
\mathsf{Hom}^{-2}_{\mathcal{B}}(X,Y)$ satisfying the relations
$d(c_X)=\mathbf{1}_X$, $d(c_Y)=\mathbf{1}_Y$ and $d(h)=c_Y\circ f +
f\circ c_X$.
Suppose now, that we dispose of a contraction $(c_1,c_2)$ of $(X,Y)$
in $\mathcal{B}\times\mathcal{B}$. We can lift this contraction by
considering $c_X=c_1$, $c_Y=c_2$ and $h=c_2\circ f \circ c_1$. This
shows that contractions lift along the dg functor $P(\mathcal{B}) \stackrel{p_0 \times
  p_1}{\longrightarrow} \mathcal{B}\times\mathcal{B}$. We dispose of the following equivalence of dg categories
$$ \mathsf{pretr}(P(\mathcal{B})) \stackrel{\sim}{\longrightarrow}
P(\mathsf{pretr}(\mathcal{B}))\,,$$
where $\mathsf{pretr}$ denotes the pre-triangulated hull of a dg category,
see~\cite{Drinfeld}. This implies that the dg functor $p_0 \times p_1$
is a fibration in the Quillen model structure described in
\cite{cras}. This proves the proposition. 
\end{proof}

Let $\mathcal{A}$ be a cofibrant dg category and $F,G:\mathcal{A}
\rightarrow \mathcal{B}$ dg functors. The dg functors $F$ and $G$ are
homotopic in the Quillen model structure described in \cite{cras} if
and only if there exists a dg functor $H:\mathcal{A} \rightarrow
P(\mathcal{B})$ that makes the following diagram commute
$$
\xymatrix{
& \mathcal{B} \\
\mathcal{A} \ar[r]^-H \ar[ur]^F \ar[dr]_G & P(\mathcal{B}) \ar[u]_{P_0}
\ar[d]^{P_1} \\
& \mathcal{B}
}
$$
, see~\cite{Hirschhornn}.

\begin{remark}\label{cic}
Remark that a dg functor $H$ as above corresponds exactly,
see~\cite{cyclic}, to:
\begin{enumerate}
\item[-] a morphism $\eta A:F(A) \rightarrow G(A)$ of
  $Z^0(\mathcal{B})$ which becomes invertible in
  $\mathsf{H}^0(\mathcal{B})$ for all $A \in \mathcal{A}$ (but which
  will not be functorial in $A$, in general) and
\item[-] a morphism of graded $k$-modules homogeneous of degree $-1$
$$ h=h(A,B): \mathsf{Hom}_{\mathcal{A}}(A,B) \rightarrow
\mathsf{Hom}_{\mathcal{B}}(F(A),G(B))\,,$$
for all $A,B \in \mathcal{A}$ such that we have
$$(\eta B)(F(f)) -(G(f)(\eta A)=d(h(f))+h(d(f))$$
and
$$h(fg)=h(f)(F(g))+(-1)^n(G(f))h(g)\,$$
for all composable morphismes $f,g$ of $\mathcal{A}$, where $f$ is of
degree $n$.
\end{enumerate}
It is shownd in~\cite{cyclic} that if we dispose of a dg functor $H$
as above and the dg category $\mathcal{B}$ is stable under cones, we
can construct a sequence of dg functors
$$ F \rightarrow I \rightarrow G[1]\,,$$
where $I(A)$ is a contractible object of $\mathcal{B}$, for all $A \in \mathcal{B}$.
\end{remark}

\section{$Q$-model structure}\label{model}

\begin{definition}
A {\em localization pair} $\mathcal{A}$ is given by a small dg category
$\mathcal{A}_1$ and a full dg subcategory $\mathcal{A}_0 \subset
\mathcal{A}_1$. A {\em morphism} $F:\mathcal{A} \rightarrow
\mathcal{B}$ of localization pairs is given by a commutative square
$$
\xymatrix{
\mathcal{A}_0 \ar[d]_{F_0} \ar@{^{(}->}[r] & \mathcal{A}_1 \ar[d]^{F_1} \\
\mathcal{B}_0 \ar@{^{(}->}[r] & \mathcal{B}_1
}
$$
of dg functors.
\end{definition}
We denote by $\mathsf{Lp}$ the category of localization pairs.

Let $\mathcal{A}$ be a localization pair.

\begin{definition}
The {\em dg quotient} of $\mathcal{A}$, see~\cite{Drinfeld}, is the dg
category $\mathcal{A}_1/\mathcal{A}_0$
obtained from $\mathcal{A}_1$ by introducing a new morphism $h_X$ of
degree $-1$ for every object $X$ of $\mathcal{A}_0$ and by imposing the
relation $d(h_X)={\bf 1}_X$.
\end{definition}

\subsection{Morita model structure}\label{secMorita}

Let $L$ be the category with two objects $0$ and $1$ and with a unique
non identity morphism $0 \rightarrow 1$.

\begin{remark}
An immediate application of Theorem $11.6.1$ from \cite{Hirschhornn}
implies that the category $\mathsf{dgcat}^L$, i.e. the category of morphisms
in $\mathsf{dgcat}$, admits a structure of cofibrantly generated model
category whose weak equivalences $W$ are the componentwise Morita
equivalences and with generating cofibrations $\mathbf{F}^L_{I}$ and
generating trivial cofibrations $\mathbf{F}^L_{J}$, where we use the
notation of \cite{Hirschhornn}:
\end{remark}
The functor $\mathbf{F}^i_?$,
$i=0,1$, from $\mathsf{dgcat}$ to $\mathsf{dgcat}^L$ is left adjoint to the evaluation functor $Ev_i$, $i=0,1$, from
$\mathsf{dgcat}^L$ to $\mathsf{dgcat}$. By definition, we have $\mathbf{F}^L_{I} =
\mathbf{F}^0_{I} \cup \mathbf{F}^1_{I}$ and  $\mathbf{F}^L_{J} =
\mathbf{F}^0_{J} \cup \mathbf{F}^1_{J}$.

The inclusion functor $U :\mathsf{Lp} \rightarrow \mathsf{dgcat}^L$
admits a left adjoint $S$ which sends an object $G:\mathcal{B}_0
\rightarrow \mathcal{B}_1$ to the localization pair formed by
$\mathcal{B}_1$ and its full dg subcategory $\mathsf{Im}\,G$.
  
\begin{proposition}\label{naive}
The category $\mathsf{Lp}$ admits a structure of cofibrantly generated model
category whose weak equivalences $W$ are the componentwise Morita
equivalences and with generating cofibrations $\mathbf{F}^L_{I}$ and
generating trivial cofibrations $\mathbf{F}^L_{J}$.
\end{proposition}

\begin{proof}
We first prove that $\mathsf{Lp}$ is complete and cocomplete.
Let $\{X_i\}_{i \in I}$ be a diagram in $\mathsf{Lp}$. We remark that 
$$ \underset{i \in I}{\mbox{colim}} \,X_i
\stackrel{\sim}{\longrightarrow} S(\underset{i \in I}{\mbox{colim}}
\,U(X_i)) \,,$$
which implies that $\mathsf{L}p$ is cocomplete. The category
$\mathsf{Lp}$ is also complete, since it is stable under products and
equalizers in $\mathsf{dgcat}^L$.
We now prove that conditions $(1)$ and $(2)$ of Theorem $11.3.2$ from \cite{Hirschhornn} are satisfied~:
\begin{enumerate}
\item[(1)] Since $S(\mathbf{F}^L_{I})=\mathbf{F}^L_{I}$ and
  $S(\mathbf{F}^L_{J})=\mathbf{F}^L_{J}$ condition $(1)$ is verified.
\item[(2)] Since the functor $U$ clearly commutes with filtered
  colimits, it is enough to prove the following: let $Y
  \stackrel{G}{\rightarrow} Z$ be an element of the set
  $\mathbf{F}^L_{J}$, $X$ an object in $\mathsf{Lp}$ and $Y \rightarrow X$ a morphism in
  $\mathsf{Lp}$. Consider the following push-out in $\mathsf{Lp}$:
$$
\xymatrix{
Y \ar[r] \ar[d]_G  \ar@{}[dr]|{\lrcorner} & X \ar[d]^{G_{\ast}} \\
Z \ar[r] & Z \underset{Y}{\coprod} X
}
$$ 
We prove that $U(G_{\ast})$ is a weak equivalence in $\mathsf{dgcat}^L$. We consider two situations:
\begin{enumerate}
\item[-] if $G$ belongs to the set $\mathbf{F}^0_{J} \subset \mathbf{F}^L_{J}$,
  then $U(G_{\ast})$ is a weak-equivalence simply because 
  $J-\mbox{cell} \subset W$ in $\mathsf{dgcat}$ , see~\cite{IMRN}\cite{erratum}\cite{cras}.
\item[-] if $G$ belongs to the set $\mathbf{F}^1_{J} \subset \mathbf{F}^L_{J}$,
  then $Ev_1(U(G_{\ast}))$ is a Morita equivalence. In particular it
  induces a quasi-isomorphism in the $\mathsf{Hom}$ spaces and since the
  $0$-component of $G_*$ is the identity on objects, the functor $Ev_0(U(G_*))$
  is also a Morita equivalence. This implies that $U(G_*)$ is a
  weak equivalence and so condition $(2)$ is proven.
\end{enumerate}
This proves the proposition.
\end{enumerate} 
\end{proof}

We will now slightly modify the previous Quillen model
structure on $\mathsf{Lp}$.

Let $\sigma$ be the morphism of localization pairs:
$$
\xymatrix{
(\mbox{End}_{\mathcal{K}}(1) \ar@{^{(}->}[r] \ar[d]_{inc}  & \mathcal{K}) \ar@{=}[d]\\
(\mathcal{K} \ar@{=}[r] & \mathcal{K})\,,
}
$$
where $\mathsf{End}_{\mathcal{K}}(1)$ is the dg algebra of
endomorphisms of the object $1$ in $\mathcal{K}$, see~\cite{cras},
and $inc$ is the natural inclusion dg-functor. Clearly $\sigma$ is a
componentwise Morita equivalence.
We write $\tilde{\mathbf{F}}^L_{I}$ resp. $\tilde{\mathbf{F}}^L_{J}$ for
the union of $\{\sigma\}$ with $\mathbf{F}^L_{I}$ resp. $\mathbf{F}^L_{J}$.

\begin{proposition}\label{sat}
The category $\mathsf{Lp}$ admits a structure of cofibrantly generated model
category whose weak equivalences $W$ are the componentwise Morita
equivalences and with generating cofibrations $\tilde{\mathbf{F}}^L_{I}$ and
generating trivial cofibrations $\tilde{\mathbf{F}}^L_{J}$.
\end{proposition}

\begin{proof}
The proof will consist in verifying that conditions
$(1) -(6)$ of Theorem $2.1.19$ from \cite{Hovey}
are satisfied. Condition $(1)$ is clear. Since the localization pair
$(\mbox{End}_{\mathcal{K}}(1) \subset \mathcal{K})$ is small in $\mathsf{Lp}$,
conditions $(2)$ and $(3)$ are also satisfied. We have 
$$ \mathbf{F}^L_{I} - \mbox{inj} = \mathbf{F}^L_{J} - \mbox{inj} \cap W$$
and so by construction
$$ \tilde{\mathbf{F}}^L_{I} - \mbox{inj} = \tilde{\mathbf{F}}^L_{J} - \mbox{inj}
\cap W\,.$$
This shows conditions $(5)$ and $(6)$. We now prove that
$\tilde{\mathbf{F}}^L_{J} - \mbox{cell} \subset W$. Since
$\mathbf{F}^L_{J}-\mbox{cell} \subset W$ it is enough to prove that
pushouts with respect to $\sigma$ belong to $W$.
Let $\mathcal{A}$ be a localization pair and
$$T: (\mbox{End}_{\mathcal{K}}(1) \subset \mathcal{K})
\rightarrow (\mathcal{A}_0 \subset \mathcal{A}_1)$$ a morphism in
$\mathsf{Lp}$. Consider the following push-out in $\mathsf{Lp}$:
$$
\xymatrix{
(\mbox{End}_{\mathcal{K}}(1) \subset \mathcal{K}) \ar[r]^-{T}
\ar[d]_{\sigma} \ar@{}[dr]|{\lrcorner} & (\mathcal{A}_0 \subset
\mathcal{A}_1) \ar[d]^R \\
(\mathcal{K} = \mathcal{K}) \ar[r] & (\mathcal{U}_0 \subset
\mathcal{U}_1) \,.
}
$$
We remark that the morphism $T$ corresponds to specifying an
object $X$ in $\mathcal{A}_0$ and a homotopy equivalence from $X$ to an
object $Y$ in $\mathcal{A}_1$. Clearly $\mathcal{U}_1 = \mathcal{A}_1$
and $\mathcal{U}_0$ identifies with the full dg-subcategory of
$\mathcal{U}_1$ whose objects are $Y$ and those of $\mathcal{A}_0$.
Since $X$ and $Y$ are homotopy equivalent, the natural dg-functor
$R_0: \mathcal{A}_0 \hookrightarrow \mathcal{U}_0$ is a
quasi-equivalence. This proves condition $(4)$.
The proposition is now proven.
\end{proof}

\begin{remark}
Remark that in this new Quillen model structure on $\mathsf{Lp}$ we dispose of
more cofibrations and less fibrations than the Quillen model structure
of proposition~\ref{naive} since the weak equivalences are the same. 
\end{remark}

From now on, by Quillen model structure on $\mathsf{Lp}$ we mean that of
proposition~\ref{sat}.

\begin{lemma}\label{fibrSat}
A localization pair $(\mathcal{A}_0 \subset \mathcal{A}_1)$ is fibrant
in $\mathsf{Lp}$ if and only if $\mathcal{A}_0$ and $\mathcal{A}_1$
are Morita fibrant dg categories and $\mathcal{A}_0$ is stable under
homotopy equivalences in $\mathcal{A}_1$.
\end{lemma}

\begin{proof}
A localization pair $(\mathcal{A}_0 \subset \mathcal{A}_1)$ is
fibrant in $\mathsf{Lp}$ if and only if for every morphism $F$ in
$\tilde{\mathbf{F}}^L_{J}$, the following extension problem in
$\mathsf{Lp}$ is solvable:
$$
\xymatrix{
X \ar[r] \ar[d]_F & (\mathcal{A}_0 \subset \mathcal{A}_1) \\
Y  \ar@{-->}[ur] & \,.
}
$$
If $F$ belongs to $\mathbf{F}^L_{J}$ this means that $\mathcal{A}_0$ and $\mathcal{A}_1$
are fibrant and if $F=\sigma$, remark that it corresponds
exactly to the
statement that $\mathcal{A}_0$ is stable under homotopy equivalences
in $\mathcal{A}_1$.
\end{proof}

\begin{lemma}\label{cofibrant}
If the localization pair $\mathcal{A}$ is
cofibrant in $\mathsf{Lp}$ then $\mathcal{A}_1$ is cofibrant in $\mathsf{dgcat}$.
\end{lemma}

\begin{proof}
We need to construct a lift to the following problem~:
$$
\xymatrix{
 & \mathcal{C} \ar@{->>}[d]^P_{\sim} \\
\mathcal{A}_1 \ar[r] & \mathcal{B} \,,
}
$$
where $P$ is a trivial fibration in $\mathsf{dgcat}$, see proposition~\ref{sat},
and $\mathcal{A}_1 \rightarrow \mathcal{B}$ is a dg-functor. Consider
the following diagram in $\mathsf{Lp}$:
$$
\xymatrix{
 & \mathbf{F}^0_{\mathcal{C}} \ar@{->>}[d]_{\sim}^{\mathbf{F}^0_P} \\
\mathcal{A} \ar[r] \ar@{.>}[ur] & \mathbf{F}^0_{\mathcal{B}}\,.
}
$$
where $\mathcal{A} \rightarrow \mathbf{F}^0_{\mathcal{B}}$ is the natural morphism of
localization pairs. Remark that $\mathbf{F}^0_P$ belongs to
$\sigma-\mbox{inj} \cap \mathbf{F}^L_{I} - \mbox{inj}$ and so is a trivial
fibration in $\mathsf{Lp}$. Since $\mathcal{A}$ is
cofibrant in $\mathsf{Lp}$ we dispose of a lifting $\mathcal{A} \rightarrow  \mathbf{F}^0_{\mathcal{C}}$ that when
restricted to the $1$-component gives us the searched lift
$\mathcal{A}_1 \rightarrow \mathcal{C}$. This proves the lemma.
\end{proof}

\subsection{$Q$-model structure}

\begin{definition}
Let $Q : \mathsf{Lp} \rightarrow \mathsf{Lp}$ be the functor that
sends a localization pair $\mathcal{A}$ to the localization pair
$$ \overline{\mathcal{A}_0} \hookrightarrow \mathcal{A}_1/\mathcal{A}_0\,,$$
 where $\overline{\mathcal{A}_0}$ is the full dg-subcategory
of $\mathcal{A}_1/\mathcal{A}_0$ whose objets are those of
$\mathcal{A}_0$.
\end{definition}

\begin{remark}
Remark that we dispose of natural morphisms
$$ \eta_{\mathcal{A}} : (\mathcal{A}_0
\subset \mathcal{A}_1) \rightarrow (\overline{\mathcal{A}_0} \subset
\mathcal{A}_1/\mathcal{A}_0)\,$$
in $\mathsf{Lp}$.
\end{remark}

\begin{definition}
A morphism of localization pairs $F: \mathcal{A} \rightarrow \mathcal{B}$ is a $Q$-{\it weak equivalence} if the induced
morphism $Q(F)$ is a weak equivalence in the Quillen model structure of
proposition~\ref{sat}. 
\end{definition}

\begin{remark}
Remark that since the objects of $\overline{\mathcal{A}_0}$
and $\overline{\mathcal{B}_0}$ are all contractible, the dg-functor
$\overline{\mathcal{A}_0} \rightarrow \overline{\mathcal{B}_0}$ is
clearly a Morita equivalence and so the morphism $F$ is a $Q$-weak
equivalence if and only if the induced dg-functor
$\mathcal{A}_1/\mathcal{A}_0 \rightarrow \mathcal{B}_1/\mathcal{B}_0$
is a Morita equivalence.
\end{remark}

\begin{definition}
A morphism in $\mathsf{Lp}$ is a {\em cofibration} if it is
  one for the Quillen model structure of Proposition~\ref{sat} and it is a $Q$-{\em fibration} if it has the right lifting property
  with respect to all cofibrations of $\mathsf{Lp}$ which are $Q$-weak equivalences.
\end{definition}

\begin{theorem}\label{main}
The category $\mathsf{Lp}$ admits a structure of Quillen model
  category whose weak equivalences are the $Q$-weak equivalences, whose
  cofibrations are the cofibrations of $\mathsf{Lp}$ and whose fibrations are
  the $Q$-fibrations.
\end{theorem}

The proof will consist in adapting the general arguments from
chapter~X from \cite{Jardine} to our situation. We start with some remarks:
\begin{enumerate}
\item[{\bf A1}] Since $k$ is a field, the conditions of theorem $3.4$
  from \cite{Drinfeld} are satisfied and so the functor $Q$ preserves weak equivalences.
\item[{\bf A2}] The morphisms of localization pairs:
$$ 
\xymatrix{
Q(\mathcal{A})
\ar@<1ex>[r]^{\eta_{Q(\mathcal{A})}}  \ar[r]_{Q(\eta_{\mathcal{A}})}  &  QQ(\mathcal{A})
}
$$
are weak equivalences in $\mathsf{Lp}$. This follows from the fact that in
both cases we are introducing contractions to objects that are already
contractible and that the functor $Q$ is the identity functor on objects.
\end{enumerate}

\begin{lemma}\label{weak fibration}
A morphism $F:\mathcal{A}
\rightarrow \mathcal{B}$ is a fibration and
a weak equivalence of $\mathsf{Lp}$ if and only if it is a $Q$-weak
equivalence and a $Q$-fibration.
\end{lemma}

\begin{proof}
Since condition {\bf A1} is verified we can use the proof of lemma
$4.3$ in chapter~X from \cite{Jardine}.
\end{proof}

\begin{nonexample}
Remark that the Quillen model structure of proposition~\ref{sat}
is not right proper, see \cite{Hirschhornn}.\newline
Let $\mathcal{B}$ be your favorite Morita fibrant dg category, whose
derived category $\mathcal{D}(\mathcal{B})$ is not trivial. In
particular the dg functor $\mathcal{B} \rightarrow 0$, where $0$
denotes the terminal object in $\mathsf{dgcat}$ is a fibration. Let $\mathcal{A}$ be the dg category with one object $1$
and whose dg algebra of endomorphisms of $1$ is $k$. Consider the
following diagram~:
$$
\xymatrix{
 & \mathcal{B} \ar@{->>}[d]^{i_0 \circ P} \\
\mathcal{A} \ar[r]_{i_{\mathcal{A}}} & 0 \coprod \mathcal{A}\,.
}
$$
Clearly $i_{\mathcal{A}}$ is a Morita equivalence and remark that the
dg functor $i_0 \circ P$ is a fibration, since the object $1$
in $\mathcal{A}$ is not contractible. This implies that in the fiber
product
$$
\xymatrix{
\emptyset \ar[d] \ar[r] \ar@{}[dr]|{\ulcorner} & \mathcal{B} \ar@{->>}[d]^{i_0 \circ P} \\
\mathcal{A} \ar[r]_{i_{\mathcal{A}}} & 0 \coprod \mathcal{A}\,,
}
$$
the dg functor $\emptyset \rightarrow \mathcal{B}$ is not a Morita
equivalence and so this Quillen model structure is not right proper.
This implies that the Quillen model structure of proposition~\ref{sat}
is also not right proper. Apply the functor $\mathbf{F}^0_?$ from
$\mathsf{dgcat}$ to $\mathsf{Lp}$ to the previous fiber product~:
$$
\xymatrix{
\emptyset = \mathbf{F}^0_{\emptyset} \ar@{}[dr]|{\ulcorner}   \ar[r] \ar[d] &
\mathbf{F}^0_{\mathcal{B}} \ar@{->>}[d]^{\mathbf{F}^0_{i_0 \circ P}} \\
\mathbf{F}^0_{\mathcal{A}} \ar[r]_{\mathbf{F}^0_{i_{\mathcal{A}}}} &
\mathbf{F}^0_{0 \coprod \mathcal{A}}\,.
}
$$
We dispose of a fiber product since the functor $\mathbf{F}^0_?$
preserves limits. Clearly $\mathbf{F}^0_{i_{\mathcal{A}}}$ is a weak
equivalence in $\mathsf{Lp}$ and remark that the morphism
$\mathbf{F}^0_{i_0\circ P}$ belongs to $\sigma-\mbox{inj} \cap
\mathbf{F}^L_{J}-\mbox{inj}$, which implies that it is a fibration in $\mathsf{Lp}$.
\end{nonexample}

Nevertheless we have the following lemma.

\begin{lemma}\label{fiber}
Let $\mathcal{A}$ be a localization pair such
that the natural morphism
$$ \eta_{\mathcal{A}}: \mathcal{A} \longrightarrow Q(\mathcal{A})$$
is a weak equivalence in $\mathsf{Lp}$. Let $F:\mathcal{W} \rightarrow
Q(\mathcal{A})$  be a fibration in $\mathsf{Lp}$. Then the
morphism
$$ \eta^*_{\mathcal{A}}: \mathcal{W} \underset{Q(\mathcal{A})}{\times} \mathcal{A} \longrightarrow \mathcal{W}$$
is a weak equivalence in $\mathsf{Lp}$.
\end{lemma}

\begin{proof}
We remark that each component of the morphism $\eta_{\mathcal{A}}$ is the identity functor on the objects of
the dg categories envolved. Since fiber products in $\mathsf{Lp}$ are
calculated componentwise, we conclude that each component of the morphism
$\eta^*_{\mathcal{A}}$ is the identity
functor on the objects. Let $X$ and $Y$ be arbitrary objects of
$\mathcal{W}_1$. We remark that we dispose of the following fiber
product in  $\mathsf{Ch}(k)$~:
$$
\xymatrix{
\mbox{Hom}_{\mathcal{W}_1
  \underset{\mathcal{A}_1/\mathcal{A}_0}{\times} \mathcal{A}_1}(X,Y)
\ar[rr] \ar[d]_{\eta^*(F_1X,F_1Y)} \ar@{}[drr]|{\ulcorner} & &
\mbox{Hom}_{\mathcal{A}_1}(F_1X,F_1Y)
\ar[d]_{\sim}^{\eta(F_1X,F_1Y)}\\
\mbox{Hom}_{\mathcal{W}_1}(X,Y) \ar@{->>}[rr]_-{F_1(X,Y)} & & \mbox{Hom}_{\mathcal{A}_1/\mathcal{A}_0}(F_1X,F_1Y)\,.
}
$$
Since $F$ is a fibration in $\mathsf{Lp}$, $F_1(X,Y)$ is a fibration in
 the projective model structure on $\mathsf{Ch}(k)$ and since this
 Quillen model structure on $\mathsf{Ch}(k)$ is right proper, $\eta^*(F_1X,F_1Y)$ is a quasi-isomorphism.
We could do the same argument for $X$ and $Y$ objects in $\mathcal{W}_0$ instead of
$\mathcal{W}_1$. This proves the lemma.
\end{proof}

\begin{lemma}\label{fibration}
Suppose that $F:\mathcal{A}_0 \rightarrow \mathcal{B}$ is a fibration in
$\mathsf{Lp}$ and that $\eta_{\mathcal{A}}$ and 
$\eta_{\mathcal{B}}$ are weak equivalences
of $\mathsf{Lp}$. Then $F$ is a $Q$-fibration.
\end{lemma}

\begin{proof}
Consider exactly the same proof as for lemma $4.4$ in chapter~X from \cite{Jardine},
but use lemma~\ref{fiber} instead of the right properness assumption
on $\mathsf{Lp}$.
\end{proof}

\begin{lemma}\label{facto}
Any morphism $F: Q(\mathcal{A})
\rightarrow Q(\mathcal{B})$ has a factorization
$F=P \circ I$ where $P:\mathcal{Z} \rightarrow Q(\mathcal{B})$ is a
$Q$-fibration and $I: Q(\mathcal{A}) \rightarrow \mathcal{Z}$ is a cofibration and a $Q$-weak equivalence.  
\end{lemma}

\begin{proof}
Since lemma~\ref{fibration} and conditions {\bf A1} and {\bf A2}
are satisfied, we consider the proof of lemma $4.5$ in chapter~X from
\cite{Jardine}.
\end{proof}

Let $\mathcal{A}$ be a localization pair. By
condition {\bf A2} we know that the natural morphism:
$$ \eta_{\mathcal{A}}:(\mathcal{A}_0 \subset
\mathcal{A}_1) \longrightarrow (\overline{\mathcal{A}_0} \subset
\mathcal{A}_1/\mathcal{A}_0)$$
is a $Q$-weak equivalence in $\mathsf{Lp}$.

\begin{lemma}\label{Q-weak}
Let $F:\mathcal{Z} \rightarrow Q(\mathcal{A})$ be a
fibration in $\mathsf{Lp}$. Then the induced morphism
$$ \eta^*_{\mathcal{A}}: \mathcal{Z} \underset{Q(\mathcal{A})} {\times} \mathcal{A} \longrightarrow \mathcal{Z}$$
is a $Q$-weak equivalence in $\mathsf{Lp}$.
\end{lemma}

\begin{proof}
We need to prove that $Q(\eta^*_{\mathcal{A}})$ is a weak equivalence in $\mathsf{Lp}$.
\begin{enumerate}
\item[(1)] We prove that the induced morphism:
$$
Q(\eta_{\mathcal{A}})^* : Q(\mathcal{Z})
\underset{QQ(\mathcal{A})}{\times} Q(\mathcal{A}) \longrightarrow Q(\mathcal{Z})
$$
is a weak equivalence in $\mathsf{Lp}$. Remark first that since $F$ is a
fibration in $\mathsf{Lp}$, the dg functors $F_0$ and $F_1$ are Morita fibrations,
see~\cite{IMRN}, and so they are surjective at the level of
$\mathsf{Hom}$-spaces. We now show that the dg functor
$F_0:\mathcal{Z}_0 \rightarrow \overline{\mathcal{A}_0}$ is surjectif
on objects. If $\overline{\mathcal{A}_0}$ is the empty dg category
then so is $\mathcal{Z}_0$ and the claim is showed. If
$\overline{\mathcal{A}_0}$ is not empty, every object $X$ in
$\overline{\mathcal{A}_0}$ is contractible and since the dg functor $F_0$ belongs to
$C-\mbox{inj}$ there exists an object $Y$ in $\mathcal{Z}_0$ such
that $F_0(Y)=X$. This implies that each component of the morphism
$$ Q(F): Q(\mathcal{Z}) \longrightarrow
QQ(\mathcal{A})$$
is a dg functor that is surjectif at the level of
$\mbox{Hom}$-spaces. Since by condition {\bf A2} the morphism
$$ (Q\eta_{\mathcal{A}}): Q(\mathcal{A}) \longrightarrow QQ(\mathcal{A})$$
is a weak equivalence an argument analogue to the proof of lemma~\ref{fiber} (we
have just proved that $F_1(X,Y)$ is a fibration in the projective
model structure on $\mathsf{Ch}(k)$),
  proves the condition $(1)$.
\item[(2)] We prove that the induced morphism:
$$ Q(\mathcal{Z}
\underset{Q(\mathcal{A})}{\times} \mathcal{A}) \longrightarrow   Q(\mathcal{Z})
\underset{QQ(\mathcal{A})}{\times} Q(\mathcal{A})$$
is an isomorphism in $\mathsf{Lp}$.
Since by construction the functor $Q$ is the identity functor on
objects, both components of the above morphism are also the identity
on objects. Let us consider de $1$-component of the above morphism. Let $X$ and
$Y$ be objects of $\mathcal{Z}_1/\mathcal{Z}_0$. We
dispose of the following fiber product in $\mathsf{Ch}(k)$~:

$$
\xymatrix{
\mbox{Hom}_{\mathcal{Z}_1/\mathcal{Z}_0
  \underset{(\mathcal{A}_1/\mathcal{A}_0)/\overline{\mathcal{A}_0}}{\times} \mathcal{A}_1/\mathcal{A}_0 }(X,Y) \ar[rr] \ar[d]  \ar@{}[drr]|{\ulcorner} & & \mbox{Hom}_{\mathcal{A}_1/\mathcal{A}_0}(F_1(X),F_1(Y)) 
\ar[d]^{Q \eta_{\mathcal{A}}} \\
\mbox{Hom}_{\mathcal{Z}_1/\mathcal{Z}_0}(X,Y) \ar@{->>}[rr]^-{QF_1} &&
\mbox{Hom}_{(\mathcal{A}_1/\mathcal{A}_0)/\overline{\mathcal{A}_0}}(F_1(X),F_1(Y))\,.
}
$$
Remark that the functor $Q \eta_{\mathcal{A}}$, resp. $QF_1$, sends the contractions in
$\mathcal{A}_1/\mathcal{A}_0$, resp. $\mathcal{Z}_1/\mathcal{Z}_0$, associated with the objects of
$\mathcal{A}_0$, resp. $\mathcal{Z}_0$, to the new contractions in
$(\mathcal{A}_1/\mathcal{A}_0)/\overline{\mathcal{A}_0}$ associated
with the objects of $\overline{\mathcal{A}_0}$.
Recall that we dispose of the following fiber product in
$\mathsf{Ch}(k)$:
$$
\xymatrix{
\mbox{Hom}_{\mathcal{Z}_1
  \underset{\mathcal{A}_1/\mathcal{A}_0}{\times} \mathcal{A}_1}(X,Y)
\ar[r] \ar[d]  \ar@{}[dr]|{\ulcorner} & \mbox{Hom}_{\mathcal{A}_1}(F_1X,F_1Y) \ar[d]^{\eta} \\
\mbox{Hom}_{\mathcal{Z}_1}(X,Y) \ar@{->>}[r]^-{F_1} &
\mbox{Hom}_{\mathcal{A}_1/\mathcal{A}_0}(F_1X, F_1Y)\,.
}
$$ 
A analysis of the above fiber products shows that the induced morphism
$$ \mbox{Hom}_{(\mathcal{Z}_1
  \underset{\mathcal{A}_1/\mathcal{A}_0}{\times}
  \mathcal{A}_1)/(\mathcal{Z}_0
    \underset{\overline{\mathcal{A}_0}}{\times}
    \mathcal{A}_0)}(X,Y) \stackrel{\sim}{\longrightarrow}
\mbox{Hom}_{\mathcal{Z}_1/\mathcal{Z}_0
  \underset{(\mathcal{A}_1/\mathcal{A}_0)/\overline{\mathcal{A}_0}}{\times} \mathcal{A}_1/\mathcal{A}_0}(X,Y)$$
is an isomorphism in $\mathsf{Ch}(k)$. The same argument applies to the
$0$-component of the above morphism. This proves condition $(2)$.
\end{enumerate}
Now, conditions $1)$ and $2)$ imply that the morphism
$$ Q(\mathcal{Z}
\underset{Q(\mathcal{A})}{\times} \mathcal{A}) \stackrel{ (Q\eta)^*_{\mathcal{A}}}{\longrightarrow} Q(\mathcal{Z})$$
is a weak equivalence in $\mathsf{Lp}$, which is exactly the statement of the
lemma. The lemma is then proved.

\end{proof}

\begin{lemma}\label{Q-fibration}
Any morphism $F:\mathcal{A}
\rightarrow \mathcal{B}$ of $\mathsf{Lp}$ has a
factorization $F=Q \circ J$ where
$Q:\mathcal{Z} \rightarrow
\mathcal{B}$ is a $Q$-fibration and
$J:\mathcal{A} \rightarrow \mathcal{Z}$ is a cofibration and a $Q$-weak equivalence.
\end{lemma}

\begin{proof}
Consider exactly the same proof as for lemma $4.6$ in chapter~X from \cite{Jardine},
but use lemma~\ref{Q-weak} instead of condition {\bf A3}.
\end{proof}

We now prove theorem~\ref{main}.
\begin{proof}
We will prove that conditions $M1-M5$ of definition $7.1.3$ from \cite{Hirschhornn}
are satisfied.
By the proof of proposition~\ref{naive}, the category $\mathsf{Lp}$ is complete and
cocomplete and so condition $M1$ is verified.
By definition the $Q$-weak equivalences in $\mathsf{Lp}$ satisfy condition
$M2$.
Clearly the $Q$-weak equivalences and $Q$-fibrations in $\mathsf{Lp}$ are
stable under retractions. Since the cofibrations are those of
proposition~\ref{sat} condition $M3$ is verified. Finally
lemma~\ref{weak fibration} implies condition $M4$ and
lemmas~\ref{weak fibration} and \ref{Q-fibration} imply condition $M5$.
\end{proof}

We denote by $\mathsf{Ho}(\mathsf{Lp})$ the homotopy category of
$\mathsf{Lp}$ given by theorem~\ref{main}.

Let $\mathcal{A}$ be a localization pair.

\begin{lemma}\label{fibrant1}
If $\mathcal{A}$ is fibrant, in the Quillen model structure of
proposition~\ref{sat}, and the morphism $\eta_{\mathcal{A}}:
\mathcal{A} \rightarrow Q(\mathcal{A})$ is a weak
equivalence in $\mathsf{Lp}$ then $\mathcal{A}$ is $Q$-fibrant.
\end{lemma}

\begin{proof}
We need to show that the morphism $\mathcal{A}
\stackrel{P}{\rightarrow} 0$ is a $Q$-fibration, where $0$ denotes the
terminal object in $\mathsf{Lp}$. Consider the
following diagram:
$$
\xymatrix{
\mathcal{A} \ar[d]_P \ar[rr]^-{\eta_{\mathcal{A}}} & & Q(\mathcal{A})
\ar[d]^{Q(P)} \\
0 \ar[rr]_-{\eta} & & Q(0) \,.
}
$$
Factorize the morphism $Q(P)$ as
$$
\xymatrix{
Q(\mathcal{A}) \ar[r]^i \ar[dr]_{Q(P)} & \mathcal{Z} \ar[d]^q\\ 
& Q(0) \,,
}
$$
where $i$ is a trivial cofibration and $q$ a fibration in
$\mathsf{Lp}$. By the proof of lemma~\ref{facto}, $q$ is a $Q$-fibration. Since the morphism $0
\rightarrow Q(0)$ is a weak equivalence, lemma~\ref{fiber} implies
that the induced morphism $0 \underset{Q(0)}{\times}
\mathcal{A} \rightarrow \mathcal{Z}$ is a weak
equivalence. Since $\eta_{\mathcal{A}}$ is a weak equivalence
the induced morphism
$$ \theta: \mathcal{A} \rightarrow 0 \underset{Q(0)}{\times}
\mathcal{Z}$$ is also a weak equivalence. Factorize the morphism $\theta$ as 
$$
\xymatrix{
\mathcal{A} \ar[r]^j \ar[dr]_{\theta} & \mathcal{W}
\ar[d]^{\pi}\\
& 0 \underset{Q(0)}{\times} \mathcal{Z}\,,
}
$$
where $\pi$ is a trivial fibration of $\mathsf{Lp}$ and $j$ is a trivial
cofibration. Then $q_* \circ \pi$ is a $Q$-fibration and the lifting
exists in the diagram~:
$$
\xymatrix{
\mathcal{A} \ar@{=}[r] \ar[d]_j & \mathcal{A} \ar[d]^p \\
\mathcal{W} \ar@{.>}[ur] \ar[r]_{q_* \circ \pi} & 0\,.
}
$$
Thus $P$ is a rectract of a $Q$-fibration, and is therefore a
$Q$-fibration itself. This proves the lemma.
\end{proof}

\begin{lemma}\label{fibrant2}
If $\mathcal{A}$ is $Q$-fibrant, then $\mathcal{A}$ is
fibrant in $\mathsf{Lp}$ and the natural morphism 
$$ \eta_{\mathcal{A}}: \mathcal{A} \rightarrow
Q(\mathcal{A})$$
is a weak equivalence.
\end{lemma}

\begin{proof}
Since the $Q$-model structure on $\mathsf{Lp}$ has less fibrations
than the Quillen model structure of proposition~\ref{sat}, the
localization pair $\mathcal{A}$ is fibrant in $\mathsf{Lp}$.
Consider the following diagram:
$$
\xymatrix{
\mathcal{A} \ar[rr]^-{\eta_{\mathcal{A}}} \ar[d]_P & & Q(\mathcal{A})
\ar[d]^{Q(P)}\\
0 \ar[rr]_{\eta} & &  Q(0) \,.
}
$$
Factorize $Q(P)=q\circ i$ as in the previous lemma. We dispose of the
following diagram~:
$$
\xymatrix{
\mathcal{A} \ar[r]^-{\theta} \ar[d]_P &  0
\underset{Q(0)}{\times} \mathcal{Z} \ar[dl]^{q_*} \\
0 &
}
$$
Since $p$ and $q_*$ are $Q$-fibrations, $\mathcal{A}$ and
$\mathcal{Z}$ are $Q$-fibrant objects in $\mathsf{Lp}$ and $\theta$
is a $Q$-weak equivalence in $\mathsf{Lp}$. By application of lemma $7.7.1$ b)
from \cite{Hirschhornn} to $\theta$ and using lemma~\ref{weak
  fibration} we conclude that $\theta$ is a weak
equivalence. Since so is $i$, we conclude that $\eta_{\mathcal{A}}$ is also a
weak equivalence. This proves the lemma.
\end{proof}

\begin{remark}
By lemmas~\ref{fibrant1} and \ref{fibrant2} a localization pair
$\mathcal{A}$ is $Q$-fibrant if and only if
it is fibrant in $\mathsf{Lp}$ and the natural morphism 
$$ \eta_{\mathcal{A}}: \mathcal{A} \longrightarrow Q(\mathcal{A})$$
is a weak equivalence. 
\end{remark}

We now describe explicitly the $Q$-fibrant objects in $\mathsf{Lp}$.

\begin{proposition}\label{Q-fibrant}
A localization pair $\mathcal{A}$ is
$Q$-fibrant, i.e. fibrant in the model structure of Theorem~\ref{main}, if and
only if it is isomorphic in $\mathsf{Lp}$ to a localization pair of the
form~:
$$ (\mathcal{B}_{contr} \subset \mathcal{B})\,,$$
where $\mathcal{B}$ is a fibrant dg category and
$\mathcal{B}_{contr}$ is the full dg subcategory of contractible
objects in $\mathcal{B}$.
\end{proposition}

\begin{proof}
Suppose first that $\mathcal{A}$ is
$Q$-fibrant. Since it is also fibrant in $\mathsf{Lp}$ the dg category
$\mathcal{A}_1$ is fibrant in $\mathsf{dgcat}$. Since the morphism
$$ \eta_{\mathcal{A}}:(\mathcal{A}_0 \subset \mathcal{A}_1) \longrightarrow
(\overline{\mathcal{A}_0} \subset \mathcal{A}_1/\mathcal{A}_0)$$
is a weak equivalence all the objects of $\mathcal{A}_0$ are
contractible. Since $\mathcal{A}$ is fibrant
in $\mathsf{Lp}$ by lemma~\ref{fibrSat} $\mathcal{A}_0$ is stable
under homotopy equivalences in $\mathcal{A}_1$. This implies that $\mathcal{A}_0$ is in fact the
full dg subcategory of contractible objects of $\mathcal{A}_1$.
Consider now a localization pair $(\mathcal{B}_{contr} \subset
\mathcal{B})$ as in the statement of the proposition.
We remark that since $\mathcal{B}$ is fibrant in $\mathsf{dgcat}$,
then $\mathcal{B}_{contr}$ it is also fibrant. Clearly $(\mathcal{B}_{contr}
\subset \mathcal{B})$ satisfies the extension condition in
what regards $\sigma$ and the morphism 
$$ \eta: (\mathcal{B}_{contr} \subset \mathcal{B}) \longrightarrow
(\overline{\mathcal{B}_{contr}} \subset
\mathcal{B}/\mathcal{B}_{contr})$$
is a weak equivalence in $\mathsf{Lp}$. This proves the proposition.
\end{proof}

\section{Closed symmetric monoidal structure}
Let $\mathcal{A}$ and $\mathcal{B}$ be small dg categories. We denote by
$\mathcal{A} \otimes \mathcal{B}$ the tensor product of $\mathcal{A}$
and $\mathcal{B}$, see \cite{dg-cat-survey} \cite{Toen}, and by
$\mathsf{Fun}_{dg}(\mathcal{A},\mathcal{B})$ the dg category of
dg functors from $\mathcal{A}$ to $\mathcal{B}$, see
\cite{dg-cat-survey} \cite{ICM}. 

\begin{definition}
The {\it internal $\mathsf{Hom}$ functor} in $\mathsf{Lp}$ 
$$ \mathsf{Hom}(-,-): \mathsf{Lp}^{op}\times \mathsf{Lp}
\longrightarrow \mathsf{Lp}\,,$$
associates to the localization pairs  $(\mathcal{A}_0 \subset \mathcal{A}_1)$, $(\mathcal{B}_0
\subset \mathcal{B}_1)$ the localization pair~:
$$ ( \mathsf{Fun}_{dg}(\mathcal{A}_1,\mathcal{B}_0) \subset
\mathsf{Fun}_{dg}(\mathcal{A}_0,\mathcal{B}_0) \underset{\mathsf{Fun}_{dg}(\mathcal{A}_0,\mathcal{B}_1)}{\times}\mathsf{Fun}_{dg}(\mathcal{A}_1,\mathsf{B}_1))\,.$$ 
\end{definition} 

\begin{definition}
The {\it tensor product} functor in $\mathsf{Lp}$
$$ - \otimes- : \mathsf{Lp} \times \mathsf{Lp} \longrightarrow
\mathsf{Lp} $$
associates to the localization pairs $(\mathcal{A}_0 \subset
\mathcal{A}_1)$, $(\mathcal{B}_0 \subset \mathcal{B}_1)$ the
localization pair~:
$$( \mathcal{A}_0\otimes\mathcal{B}_1 \cup \mathcal{A}_1\otimes
\mathcal{B}_0 \subset \mathcal{A}_1\otimes \mathcal{B}_1)\,,$$
where  $\mathcal{A}_0\otimes\mathcal{B}_1 \cup \mathcal{A}_1\otimes
\mathcal{B}_0$ is the full dg subcategory of $\mathcal{A}_1\otimes
\mathcal{B}_1$ consisting of those objects $a\otimes b$ of
$\mathcal{A}_1 \otimes \mathcal{B}_1$ such that $a$ belongs to
$\mathcal{A}_0$ or $b$ belongs to $\mathcal{B}_0$.
\end{definition}

Let $\mathcal{A}=(\mathcal{A}_0 \subset \mathcal{A}_1)$,
$\mathcal{B}=(\mathcal{B}_0 \subset \mathcal{B}_1)$ and
$\mathcal{C}=(\mathcal{C}_0 \subset \mathcal{C}_1)$ be
localization pairs.

\begin{proposition}
The category $\mathsf{Lp}$ endowed with the functors
$\mathsf{Hom}(-,-)$ and $-\otimes-$ is a closed symmetric monoidal
category. In particular we dispose of a natural isomorphism in $\mathsf{Lp}$:
$$ \mathsf{Hom}_{\mathsf{Lp}}(\mathcal{A} \otimes
\mathcal{B}, \mathcal{C}) \stackrel{\sim}{\longrightarrow}
\mathsf{Hom}_{\mathsf{Lp}}(\mathcal{A},
\mathsf{Hom}(\mathcal{B}, \mathcal{C}))\,.$$
\end{proposition}

\begin{proof}
Consider the following commutative square in $\mathsf{dgcat}$~:
$$
\xymatrix{
\mathcal{A}_0 \ar[d] \ar[rr] & & \mathsf{Fun}_{dg}(\mathcal{B}_1,
\mathcal{C}_0) \ar[d] \\
\mathcal{A}_1 \ar[rr] & &
\mathsf{Fun}_{dg}(\mathcal{B}_0, \mathcal{C}_0)
\underset{\mathsf{Fun}_{dg}(\mathcal{B}_0,\mathcal{C}_1)}{\times}
\mathsf{Fun}_{dg}(\mathcal{B}_1,\mathcal{C}_1)\,,
}
$$
which corresponds exactly to an element of
$\mathsf{Hom}_{\mathsf{Lp}}(\mathcal{A},\mathsf{Hom}(\mathcal{B},\mathcal{C}))$.
Recall from \cite{dg-cat-survey} that $\mathsf{dgcat}$ endowed with
  the functors $-\otimes-$ and $\mathsf{Hom}(-,-)$ is a closed
  symmetric monoidal category. This implies by adjunction that the
  commutative square above corresponds to the following commutative
  square in $\mathsf{dgcat}$:
$$
\xymatrix{
\mathcal{A}_0\otimes\mathcal{B}_1 \underset{\mathcal{A}_0 \otimes
  \mathcal{B}_0}{\times} \mathcal{A}_1\otimes \mathcal{B}_0 \ar[d]
\ar[r] & \mathcal{C}_0 \ar[d] \\
\mathcal{A}_1\otimes \mathcal{B}_1 \ar[r] & \mathcal{C}_1\,.
}
$$
This commutative square can be seen simply as a morphism in
$\mathsf{dgcat}^L$ from 
$$\mathcal{A}_0\otimes\mathcal{B}_1 \underset{\mathcal{A}_0 \otimes
  \mathcal{B}_0}{\times} \mathcal{A}_1\otimes \mathcal{B}_0
\longrightarrow \mathcal{A}_1 \otimes \mathcal{B}_1$$
 to the localization pair $(\mathcal{C}_0
\subset \mathcal{C}_1)$. Remark that the morphism 
$$\mathcal{A}_0\otimes\mathcal{B}_1 \underset{\mathcal{A}_0 \otimes
  \mathcal{B}_0}{\times} \mathcal{A}_1\otimes \mathcal{B}_0
\rightarrow \mathcal{A}_1 \otimes \mathcal{B}_1$$ 
of dg categories is
injective on objects and that its image consists of those objects
$a\otimes b$ of $\mathcal{A}_1 \otimes \mathcal{B}_1$ such that $a$
belongs to $\mathcal{A}_0$ or $b$ belongs to $\mathcal{B}_0$.
This implies that 
$$ \mathsf{Im}\,(\mathcal{A}_0\otimes\mathcal{B}_1 \underset{\mathcal{A}_0 \otimes
  \mathcal{B}_0}{\times} \mathcal{A}_1\otimes \mathcal{B}_0
\rightarrow \mathcal{A}_1 \otimes \mathcal{B}_1) = \mathcal{A}
\otimes \mathcal{B}\,,$$
and by the adjunction $(S,U)$ from subsection~\ref{secMorita}, this
last commutative square in $\mathsf{dgcat}$ corresponds exactly to an
element of
$\mathsf{Hom}_{\mathsf{Lp}}(\mathcal{A}\otimes\mathcal{B},
\mathcal{C})$. This proves the proposition.
\end{proof}

\begin{remark}
Remark that the unit object is the localization pair $(\emptyset
\subset \mathcal{A})$, where $\mathcal{A}$ is the dg category with 
one object and whose dg algebra of endomorphisms is $k$.
\end{remark}

\section{Derived internal \mbox{Hom}-functor}

Let $\mathcal{A}$ be a cofibrant dg category and $\lambda$ an infinite
cardinal whose size is greater or equal to the cardinality of the set of isomorphism classes of objects in the category
$\mathsf{H}^0(\mathcal{A})$.
Let $\mathcal{B}$ be a Morita fibrant dg category.

\begin{definition}\label{lambda}
Let $\mathcal{B}_{\lambda}$ be the full dg subcategory of
$\mathcal{C}_{dg}(\mathcal{B})$, whose objects are:
\begin{itemize}
\item[-] the right $\mathcal{B}$ dg modules $M$ such that $M \oplus D$
  is representable for a contractible right $\mathcal{B}$ dg module
  $D$ and
\item[-] the right $\mathcal{B}$ dg modules of the form $\widehat{B}\oplus C$, where $B$
is an object of $\mathcal{B}$ and the right $\mathcal{B}$ dg module
$C$ is a direct factor of $\underset{i \in I}{\bigoplus}
\mathsf{cone}(\mathbf{1}_{\widehat{B_i}})$, with $B_i$ an object of
$\mathcal{B}$ and $I$ a set of cardinality bounded by $\lambda$.
\end{itemize}
\end{definition}

Let $\mathsf{rep}_{dg}(\mathcal{A},\mathcal{B})$ be the dg category as in
\cite{dg-cat-survey} \cite{Toen}.

\begin{remark}
Remark that we dispose of a quasi-equivalence $\mathcal{B}
\stackrel{h}{\rightarrow} \mathcal{B}_{\lambda}$ and that the objects
  of $\mathcal{B}_{\lambda}$ are cofibrant and quasi-representable as
  right dg $\mathcal{B}$ modules, see~\cite{Toen}. This implies
  that we dispose of a natural dg functor:
$$ 
\overline{\mathsf{Fun}}_{dg}(\mathcal{A}, \mathcal{B}_{\lambda}) :=
\mathsf{Fun}_{dg}(\mathcal{A}, \mathcal{B}_{\lambda}) /
\mathsf{Fun}_{dg}(\mathcal{A}, (\mathcal{B}_{\lambda})_{contr}) \stackrel{\Phi}{\longrightarrow}
\mathsf{rep}_{dg}(\mathcal{A},\mathcal{B})\,.
$$
\end{remark}

\begin{theorem}\label{rep1}
The natural induced dg functor:
$$
\mathsf{Fun}_{dg}(\mathcal{A},\mathcal{B}_{\lambda})/\mathsf{Fun}_{dg}(\mathcal{A},(\mathcal{B}_{\lambda})_{contr})
\stackrel{\Phi}{\longrightarrow}
\mathsf{rep}_{dg}(\mathcal{A},\mathcal{B})\,,$$
is a quasi-equivalence. 
\end{theorem}

\begin{proof}
We prove first that $\mathsf{H}^0(\Phi)$ is essentially surjective. We
dispose of the following composition of dg functors
$$ 
\mathsf{Fun}_{dg}(\mathcal{A},\mathcal{B}) \stackrel{I}{\longrightarrow}
\overline{\mathsf{Fun}}_{dg}(\mathcal{A}, \mathcal{B}_{\lambda})
\stackrel{\Phi}{\longrightarrow}
\mathsf{rep}_{dg}(\mathcal{A},\mathcal{B})\,.
$$
Since $\mathcal{A}$ is a cofibrant dg category, lemma $4.3$ and sub-lemma $4.4$
from \cite{Toen} imply that $\mathsf{H}^0(\Phi \circ I)$ is essentially
surjective and so we conclude that so is $\mathsf{H}^0(\Phi)$.

We now prove also that the functor $\mathsf{H}^0(I)$ is essentially surjective.
Let $F:\mathcal{A} \rightarrow \mathcal{B}_{\lambda}$ be a dg
functor. Since $\mathcal{A}$ is a cofibrant dg category and $h$ is a
quasi-equivalence, there exists a dg functor $F':\mathcal{A}
\rightarrow \mathcal{B}$ such that $F$ and $h\circ F'$ are homotopic
in the Quillen model structure constructed in~\cite{cras}. Remark that
since $\mathcal{B}$ is a Morita fibrant dg category so is
$\mathcal{B}_{\lambda}$. In particular $\mathcal{B}_{\lambda}$ is
stable under cones up to homotopy,
see~\cite{IMRN}\cite{erratum}. Since a cone can be obtained from a
cone up to homotopy, by adding or factoring out contractible modules, we
conclude that by definition, $\mathcal{B}_{\lambda}$ is also stable
under cones. By remark~\ref{cic} we dipose of a sequence of dg functors
$$ F \longrightarrow I \longrightarrow h \circ F'[1]\,,$$
such that $I$ belongs to
$\mathsf{Fun}_{dg}(\mathcal{A},(\mathcal{B}_{\lambda})_{contr})$. This
implies that $F$ and $h \circ F'$ become isomorphic in
$\mathsf{H}^0(\overline{\mathsf{Fun}}_{dg}(\mathcal{A},
\mathcal{B}_{\lambda}))$. This proves that the functor
$\mathsf{H}^0(\Phi)$ is essentially surjective.
Let us now prove that the functor $\mathsf{H}^0(\Phi)$ is fully faithful.
Let $F$ belong to
$\mathsf{Fun}_{dg}(\mathcal{A},\mathcal{B}_{\lambda})$. Since
$\mathsf{H}^0(I)$ is essentially surjective, we can consider $F$ as
belonging to $\mathsf{Fun}_{dg}(\mathcal{A},\mathcal{B})$. We will construct a morphism of dg functors
$$ F' \stackrel{\mu}{\longrightarrow} F\,,$$
where $\mu$ becomes invertible in $
\mathsf{H}^0(\overline{\mathsf{Fun}}_{dg}(\mathcal{A}, \mathcal{B}_{\lambda}))$
and $F'$ belongs to the left-orthogonal of the category
$\mathsf{H}^0(\mathsf{Fun}_{dg}(\mathcal{A},(\mathcal{B}_{\lambda})_{contr}))$.
Consider the $\mathcal{A}$-$\mathcal{B}$-bimodule $X_F$ naturally
associated to $F$. Consider $X_F$ as a left $\mathcal{A}$-module and
let $\mathbf{P}X_F$ denote the bar resolution of $X_F$. Remark that
$\mathbf{P}X_F$ is naturally a right $\mathcal{B}$-module and that it
is cofibrant in the projective model structure on the category of $\mathcal{A}$-$\mathcal{B}$-bimodules. Let $A$ be an object of
$\mathcal{A}$. Since the dg category $\mathcal{A}$ is cofibrant in
$\mathsf{dgcat}$, $(\mathbf{P}X_F)(?,A)$ is cofibrant as a
$\mathcal{B}$-module. We dispose of the following homotopy
equivalence
$$
\xymatrix{
(\mathbf{P}X_F)(?,A) \ar@{->>}[r]_{\sim}^{\mu_A} & X_F(?,A)\,,
}
$$
since both $\mathcal{B}$-modules are cofibrant.
This implies that the $\mathcal{B}$-module $(\mathbf{P}X_F)(?,A)$ is
isomorphic to a direct sum $X_F(?,A) \oplus C$, where $C$ is a
contractible and cofibrant $\mathcal{B}$-module.
The $\mathcal{B}$-module $C$ is in fact isomorphic to a direct factor
of a $\mathcal{B}$-module
$$ \underset{i \in I}\bigoplus(\mathsf{cone}\mathbf{1}_{\widehat{B_i}})[n_i],$$
where $I$ is a set whose cardinality is bounded by $\lambda$, $B_i$, $i \in I$ is an object of $\mathcal{B}$ and
$n_i$, $i \in I$ is an integer, see \cite{dg-cat}.

This implies, by definition of $\mathcal{B}_{\lambda}$, that the $\mathcal{B}$-module
$$ X_F(?,A)\oplus C$$ belongs to $\mathcal{B}_{\lambda}$ and so the $\mathcal{A}$-$\mathcal{B}$-bimodule $\mathbf{P}X_F$ is in fact
isomorphic to $X_{F'}$ for a dg functor $F': \mathcal{A} \rightarrow \mathcal{B}_{\lambda}$. Remark that the previous construction is
functorial in $A$ and so we dispose of a morphism of dg functors
$$ F' \stackrel{\mu}{\longrightarrow} F\,.$$
Since for each $A$ in $\mathcal{A}$, the morphism $\mu_A: F'A
\rightarrow FA$ is a retraction with contractible kernel, the morphism
$\mu$ becomes invertible in
$$\mathsf{H}^0(\overline{\mathsf{Fun}}_{dg}(\mathcal{A}, \mathcal{B}_{\lambda}))\,.$$
Let now $G$ belong to
$\mathsf{Fun}_{dg}(\mathcal{A},(\mathcal{B}_{\lambda})_{contr})$. We remark that 
$$
\mathsf{Hom}_{\mathsf{H}^0(\mathsf{Fun}_{dg}(\mathcal{A},\mathcal{B}_{\lambda}))}(F',G)
\stackrel{\sim}{\longrightarrow}
\mathsf{Hom}_{\mathcal{H}(\mathcal{A}^{op}\otimes
  \mathcal{B})}(\mathbf{P}X_F,X_G)\,,$$
where $\mathcal{H}(\mathcal{A}^{op}\otimes
  \mathcal{B})$ denotes the homotopy category of
  $\mathcal{A}$-$\mathcal{B}$ bimodules. 
Since $\mathbf{P}X_F$ is a cofibrant $\mathcal{A}$-$\mathcal{B}$-bimodule
and $X_G(?,A)$ is a contractible $\mathcal{B}$-module, for every
object $A$ in $\mathcal{A}$, we conclude that the right hand side
vanishes and $F'$ belongs to the
left-orthogonal of
$\mathsf{H}^0(\mathsf{Fun}_{dg}(\mathcal{A},(\mathcal{B}_{\lambda})_{contr}))$.
This implies that the induced functor 
$$
\mathsf{H}^0(\mathsf{Fun}_{dg}(\mathcal{A},\mathcal{B}_{\lambda})/\mathsf{Fun}_{dg}(\mathcal{A},(\mathcal{B}_{\lambda})_{contr}))
\rightarrow \mathsf{H}^0(\mathsf{rep}_{dg}(\mathcal{A},\mathcal{B}))$$
is fully faithful.
This proves the theorem.
\end{proof}

\begin{proposition}\label{rep}
The internal $\mathsf{Hom}$ functor
$$ \mathsf{Hom}(-,-): \mathsf{Lp}^{op} \times \mathsf{Lp} \rightarrow
  \mathsf{Lp} \,,$$
admits a total right derived functor
$$ \mathcal{R}\mathsf{Hom}(-,-):\mathsf{Ho}(\mathsf{Lp})^{op} \times
\mathsf{Ho}(\mathsf{Lp}) \rightarrow \mathsf{Ho}(\mathsf{Lp})\,
$$
as in definition $8.4.7$ from~\cite{Hirschhornn}.
\end{proposition}

\begin{proof}
Let $\mathcal{A}$ and $\mathcal{B}$ be localization pairs. We are now
going to define $\mathcal{R}\mathsf{Hom}(\mathcal{A},\mathcal{B})$ and
the morphism $\epsilon$ as in definition $8.4.7$
from~\cite{Hirschhornn}. We denote by
$\mathcal{A}_c \stackrel{P}{\rightarrow} \mathcal{A}$ a functorial
cofibrant resolution of $\mathcal{A}$ in $\mathsf{Lp}$ and  by
$\mathcal{B} \stackrel{I}{\rightarrow} \mathcal{B}_f$ a functorial
$Q$-fibrant resolution of $\mathcal{B}$ in $\mathsf{Lp}$.
Remember, that by proposition~\ref{Q-fibrant}, $\mathcal{B}_f$ is of
the form
$$\mathcal{B}_f=((\mathcal{B}_f)_{contr} \subset \mathcal{B}_f)\,,$$
where $\mathcal{B}_f$ is a Morita fibrant dg category.
Let $\lambda$ be an infinite cardinal whose size is greater or equal to
the cardinality of the set of isomorphim classes in the category
$\mathsf{H}^0((\mathcal{A}_c)_1)$.
Consider now the following localization pair
$$ (\mathcal{B}_f)_{\lambda} := (((\mathcal{B}_f)_{\lambda})_{contr}
\subset (\mathcal{B}_f)_{\lambda})\,,$$
where $(\mathcal{B}_f)_{\lambda}$ is as in definition~\ref{lambda}. Remark that we
dispose of a canonical weak equivalence in $\mathsf{Lp}$
$$ \mathcal{B}_f \stackrel{F}{\longrightarrow}
(\mathcal{B}_f)_{\lambda}\,.$$
We now define $\mathcal{R}\mathsf{Hom}(\mathcal{A},\mathcal{B})$ as
$\mathsf{Hom}(\mathcal{A}_c, (\mathcal{B}_f)_{\lambda})$ and we consider
for morphism $\epsilon$ the image in $\mathsf{H}^0(\mathsf{Lp})$ of the following $Q$-equivalence in
$\mathsf{Lp}$
$$ \eta : (\mathcal{A},\mathcal{B}) \stackrel{(P,I)}{\longrightarrow}
(\mathcal{A}_c, \mathcal{B}_f) \stackrel{(Id,F)}{\longrightarrow}
(\mathcal{A}_c, (\mathcal{B}_f)_{\lambda})$$
under the functor $\mathsf{Hom}(-,-)$.
We will now show that the dg category associated with the localization
pair $\mathcal{R}\mathsf{Hom}(\mathcal{A},\mathcal{B})$ is canonically
Morita equivalent to
$$ \mathsf{rep}_{dg}((\mathcal{A}_c)_1/(\mathcal{A}_c)_0
,\mathcal{B}_f)\,.$$
Remark that since $\mathcal{A}_c$ is a cofibrant object in
$\mathsf{Lp}$, by lemma~\ref{cofibrant}, $(\mathcal{A}_c)_1$ is
cofibrant in $\mathsf{dgcat}$ and so  we dispose of an exact sequence
in $\mathsf{Hmo}$, see~\cite{dg-cat-survey}
$$ (\mathcal{A}_c)_0 \hookrightarrow (\mathcal{A}_c)_1 \rightarrow
(\mathcal{A}_c)_1 / (\mathcal{A}_c)_0\,.$$
Since the dg category $(\mathcal{B}_f)$ is Morita fibrant,
the application of the functor
$\mathsf{rep}_{dg}(-,\mathcal{B}_f)$ to the previous exact
sequence induces a new exact sequence in $\mathsf{Hmo}$
$$ \mathsf{rep}_{dg}((\mathcal{A}_c)_0, \mathcal{B}_f)
\leftarrow \mathsf{rep}_{dg}((\mathcal{A}_c)_1,
\mathcal{B}_f) \leftarrow
\mathsf{rep}_{dg}((\mathcal{A}_c)_1/ (\mathcal{A}_c)_0,
\mathcal{B}_f)\,.
$$
Remember that:
$$\mathsf{Hom}(\mathcal{A}_c, (\mathcal{B}_f)_{\lambda})_1 = \mathsf{Fun}_{dg}((\mathcal{A}_c)_0,((\mathcal{B}_f)_{\lambda})_{contr})
\underset{\mathsf{Fun}_{dg}((\mathcal{A}_c)_0,(\mathcal{B}_f)_{\lambda})}{\times}
\mathsf{Fun}_{dg}((\mathcal{A}_c)_1,(\mathcal{B}_f)_{\lambda})\,.$$
Now, since the dg categories $(\mathcal{A}_c)_1$ and
$(\mathcal{B}_f)_{\lambda}$ satisfy the conditions of
theorem~\ref{rep1}, we dispose of a natural inclusion of dg categories
$$
\mathsf{Hom}(\mathcal{A}_c, (\mathcal{B}_f)_{\lambda})_1 /\mathsf{Fun}_{dg}((\mathcal{A}_c)_1,((\mathcal{B}_f)_{\lambda})_{contr})
\longrightarrow \mathsf{rep}_{dg}((\mathcal{A}_c)_1 ,\mathcal{B}_f)\,.
$$
Now remark that this inclusion induces the following Morita
equivalence
$$
\mathsf{Hom}(\mathcal{A}_c, (\mathcal{B}_f)_{\lambda})_1 /\mathsf{Fun}_{dg}((\mathcal{A}_c)_1,((\mathcal{B}_f)_{\lambda})_{contr})
\stackrel{\sim}{\longrightarrow} \mathsf{rep}_{dg}((\mathcal{A}_c)_1/(\mathcal{A}_c)_0 ,\mathcal{B}_f)\,.
$$
We now show that the functor $\mathcal{R} \mathsf{Hom}(-,-)$ preserves
$Q$-weak equivalences in $\mathsf{Lp}^{op}\times
\mathsf{Lp}$. Consider a $Q$-weak equivalence
$$ (\mathcal{A},\mathcal{B}) \rightarrow (\tilde{\mathcal{A}},\tilde{\mathcal{B}})\,,$$
in $\mathsf{Lp}^{op} \times \mathsf{Lp}$. By construction it will
induce a Morita equivalence
$$ (\tilde{\mathcal{A}}_c)_1 / (\tilde{\mathcal{A}}_c)_0
\stackrel{\sim}{\longrightarrow} (\mathcal{A}_c)_1 /
(\mathcal{A}_c)_0$$
and also a Morita equivalence
$$ \mathcal{B}_f \stackrel{\sim}{\longrightarrow}
\tilde{\mathcal{B}}_f\,.$$
This implies that the induced dg functor
$$ \mathsf{rep}_{dg}((\mathcal{A}_c)_1 /(\mathcal{A}_c)_0,
\mathcal{B}_f) \stackrel{\sim}{\longrightarrow} \mathsf{rep}_{dg}((\tilde{\mathcal{A}}_c)_1 /(\tilde{\mathcal{A}}_c)_0,\tilde{\mathcal{B}}_f)$$
is a Morita equivalence. Now, remark that we dispose of the natural zig-zag of $Q$-weak equivalences in
$\mathsf{Lp}$:
\begin{eqnarray*}
& (\mathsf{Fun}_{dg}((\mathcal{A}_c)_1,
((\mathcal{B}_f)_{\lambda})_{contr}) \subset
\mathsf{Hom}(\mathcal{A}_c, (\mathcal{B}_f)_{\lambda})_1 & \\
& \downarrow & \\
& (\overline{\mathsf{Fun}_{dg}((\mathcal{A}_c)_1,
((\mathcal{B}_f)_{\lambda})_{contr})} \subset  \mathsf{Hom}(\mathcal{A}_c, (\mathcal{B}_f)_{\lambda})_1 / \mathsf{Fun}_{dg}((\mathcal{A}_c)_1,
((\mathcal{B}_f)_{\lambda})_{contr})) & \\
& \uparrow & \\
& (\emptyset \subset \mathsf{Hom}(\mathcal{A}_c, (\mathcal{B}_f)_{\lambda})_1 / \mathsf{Fun}_{dg}((\mathcal{A}_c)_1,
((\mathcal{B}_f)_{\lambda})_{contr})) &
\end{eqnarray*}
This allow us to conclude that the  the functor $\mathcal{R} \mathsf{Hom}(-,-)$ preserves
$Q$-weak equivalences in $\mathsf{Lp}^{op}\times
\mathsf{Lp}$.
This proves the proposition.
\end{proof}

\begin{lemma}\label{tensorprod}
Let $\mathcal{A}$ be a cofibrant object in $\mathsf{Lp}$. The induced
internal tensor product functor
$$ \mathcal{A}\otimes- : \mathsf{Lp} \longrightarrow \mathsf{Lp}\,,$$
preserves $Q$-weak equivalences.
\end{lemma}

\begin{proof}
Let $F:\mathcal{B} \rightarrow \mathcal{C}$ be a $Q$-weak equivalence
in $\mathsf{Lp}$ between cofibrant objects. We prove that the induced
morphism in $\mathsf{Lp}$ 
$$ \mathcal{A} \otimes \mathcal{B} \stackrel{F_*}{\longrightarrow}
\mathcal{A}\otimes \mathcal{C}\,,$$
is a $Q$-weak equivalence. By lemma~\ref{cofibrant}, $\mathcal{A}_1$,
$\mathcal{B}_1$ and $\mathcal{C}_1$ are cofibrant dg categories in
$\mathsf{dgcat}$ and so we dispose of a morphism of exact sequences in
$\mathsf{Hmo}$:
$$
\xymatrix{
\mathcal{B}_0 \ar[d] \ar@{^{(}->}[r] & \mathcal{B}_1 \ar[d] \ar[r] &
\mathcal{B}_1/\mathcal{B}_0 \ar[d]^{\sim} \\
\mathcal{C}_0 \ar@{^{(}->}[r] & \mathcal{C}_1 \ar[r] & \mathcal{C}_1/\mathcal{C}_0\,,
}
$$
where the last column is a Morita equivalence. Since $\mathcal{A}_1$
is cofibrant in $\mathsf{dgcat}$, by applying the functor
$\mathcal{A}\otimes-$, we obtain the following morphism of exact
sequences in $\mathsf{Hmo}$:
$$
\xymatrix{
\mathcal{A}_1 \otimes \mathcal{B}_0 \ar[d] \ar[r] &
\mathcal{A}_1\otimes \mathcal{B}_1 \ar[d] \ar[r] &
\mathcal{A}_1 \otimes \mathcal{B}_1/\mathcal{B}_0 \ar[d]^{\sim} \\
\mathcal{A}_1 \otimes \mathcal{C}_0 \ar[r] & \mathcal{A}_1
\otimes \mathcal{C}_1 \ar[r] & \mathcal{A}_1 \otimes \mathcal{C}_1/\mathcal{C}_0\,.
}
$$
This implies that we dispose of the following Morita equivalence:
$$ (\mathcal{A}_1 \otimes \mathcal{B}_1)/(\mathcal{A}_1 \otimes
\mathcal{B}_0) \stackrel{\sim}{\longrightarrow} (\mathcal{A}_1 \otimes \mathcal{C}_1)/(\mathcal{A}_1 \otimes
\mathcal{C}_0)\,.$$
Let $\mathcal{H}$ be the full dg subcategory of $(\mathcal{A}_1 \otimes \mathcal{B}_1)/(\mathcal{A}_1 \otimes
\mathcal{B}_0)$, whose objects are $a \otimes b$, where $a$ belongs to
$\mathcal{A}_0$ and $\mathcal{P}$ the full dg subcategory of $(\mathcal{A}_1 \otimes \mathcal{C}_1)/(\mathcal{A}_1 \otimes
\mathcal{C}_0)$ whose objects are $a \otimes c$, where $a$ belongs to
$\mathcal{A}_0$.
We dispose of the following diagram:
$$
\xymatrix{
\mathcal{H} \ar[d]^{\sim} \ar@{^{(}->}[r] &  (\mathcal{A}_1 \otimes \mathcal{B}_1)/(\mathcal{A}_1 \otimes
\mathcal{B}_0) \ar[d]^{\sim}\\
\mathcal{P} \ar@{^{(}->}[r] & (\mathcal{A}_1 \otimes \mathcal{C}_1)/(\mathcal{A}_1 \otimes
\mathcal{C}_0)\,.
}
$$
Remark that the dg categories $\mathcal{A}\otimes \mathcal{B}$ and
$\mathcal{A}\otimes \mathcal{C}$ are Morita equivalent dg
subcategories of $((\mathcal{A}_1 \otimes \mathcal{B}_1)/(\mathcal{A}_1 \otimes
\mathcal{B}_0))/\mathcal{H}$, resp. $((\mathcal{A}_1 \otimes \mathcal{C}_1)/(\mathcal{A}_1 \otimes
\mathcal{C}_0))/\mathcal{P}$ and so we have the commutative square:
$$
\xymatrix{
((\mathcal{A}_1 \otimes \mathcal{B}_1)/(\mathcal{A}_1 \otimes
\mathcal{B}_0))/\mathcal{H} \ar[d]^{\sim} &
\mathcal{A} \otimes \mathcal{B} \ar[l]^-{\sim} \ar[d]_{\sim}^{F^*}\\
((\mathcal{A}_1 \otimes \mathcal{C}_1)/(\mathcal{A}_1 \otimes
\mathcal{C}_0))/\mathcal{P} & \mathcal{A}\otimes \mathcal{C} \ar[l]^-{\sim}\,.
}
$$
This implies the lemma.
\end{proof}

\begin{remark}
Since the internal tensor product $-\otimes-$ is symmetric,
lemma~\ref{tensorprod} implies that the total left derived functor $-\otimes-$
$$ - \overset{\mathbb{L}}{\otimes}- : \mathsf{Ho}(\mathsf{Lp})\times
\mathsf{Ho}(\mathsf{Lp}) \rightarrow \mathsf{Ho}(\mathsf{Lp})$$
exists, see definition $8.4.7$ of \cite{Hirschhornn}.
\end{remark}

\section{Relation with $\mathsf{dgcat}$}

We dispose of the following adjunction:
$$
\xymatrix{
\mathsf{Lp} \ar@<1ex>[d]^{Ev_1}\\
\mathsf{dgcat} \ar@<1ex>[u]^F \,,
}
$$
where $Ev_1$ is the evaluation functor on the $1$-component and $F$
associates to a dg category $\mathcal{A}$ the localization pair $(\emptyset
\subset \mathcal{A})$.
\begin{lemma}
If we consider on $\mathsf{dgcat}$ the Quillen model structure of~\cite{IMRN}\cite{erratum} and on $\mathsf{Lp}$ the Quillen model structure of theorem~\ref{main}, the previous
adjunction is a Quillen equivalence, see~\cite{Hirschhornn}.
\end{lemma}

\begin{proof}
The functor $F$ clearly sends Morita equivalences to weak
equivalences. By lemma~\ref{fiber} the evaluation functor $Ev_1$
preserves trivial fibrations. This shows that $F$ is a left Quillen functor.
Let $\mathcal{A}$ be a cofibrant object in $\mathsf{dgcat}$ and
$(\mathcal{B}_{contr} \subset \mathcal{B})$ a $Q$-fibrant object in
$\mathsf{Lp}$. Let $\mathcal{A} \stackrel{F}{\rightarrow} \mathcal{B}$
be a dg-functor in $\mathsf{dgcat}$. We need to show that $F$ is a
Morita equivalence if and only if the induced morphism of localization
pairs $(\emptyset \subset \mathcal{A} ) \rightarrow
(\mathcal{B}_{contr} \subset \mathcal{B})$ is a $Q$-weak
equivalence. But since the dg functor $\mathcal{B} \rightarrow
\mathcal{B}/\mathcal{B}_{contr}$ is a Morita equivalence this
automatically follows.
\end{proof}

\begin{lemma}\label{Quilleneq}
The total derived functors, $-\overset{\mathbb{L}}{\otimes}-$ and
$\mathcal{R}\mathsf{Hom}(-,-)$ in the category
$\mathsf{Ho}(\mathsf{Lp})$ correspond, under the equivalence:
$$
\xymatrix{
\mathsf{Ho}(\mathsf{Lp}) \ar@<1ex>[d]^{\mathcal{R}Ev_1} \\
\mathsf{Ho}(\mathsf{dgcat}) \ar@<1ex>[u]^F
}
$$  
to the functors, $-\overset{\mathbb{L}}{\otimes}-$  and
$\mathsf{rep}_{dg}(-,-)$, see \cite{dg-cat-survey}\cite{Toen}, in the
category $\mathsf{Ho}(\mathsf{dgcat})$.
\end{lemma}

\begin{proof}
Let $\mathcal{A}$ and $\mathcal{B}$ be dg categories. Then
$\mathcal{A} \overset{\mathbb{L}}{\otimes}\mathcal{B}$ identifies with
$\mathcal{A}_c \otimes \mathcal{B}$, where $\mathcal{A}_c$ is a
cofibrant resolution of $\mathcal{A}$ in $\mathsf{dgcat}$. Since
$F(\mathcal{A}_c)$ is cofibrant in $\mathsf{Lp}$ by lemma~\ref{Quilleneq}, we
have the following zig-zag:
$$F(\mathcal{A}) \overset{\mathbb{L}}{\otimes} F(\mathcal{B})
\stackrel{\sim}{\leftarrow} F(\mathcal{A}_c) \otimes F(\mathcal{B})=
F(\mathcal{A}_c \overset{\mathbb{L}}{\otimes} \mathcal{B})
\stackrel{\sim}{\rightarrow} F(\mathcal{A}\otimes \mathcal{B})\,,$$
of weak equivalences in $\mathsf{Lp}$. This proves that the total left derived tensor
products in $\mathsf{Ho}(\mathsf{Lp})$ and
$\mathsf{Ho}(\mathsf{dgcat})$ are identified.
Now, $\mathsf{rep}_{dg}(\mathcal{A},\mathcal{B})$ identifies with
$\mathsf{rep}_{dg}(\mathcal{A}_c,\mathcal{B}_f)$, where
$\mathcal{B}_f$ is a fibrant resolution of
$\mathcal{B}$ in $\mathsf{dgcat}$. By definition
$$ \mathcal{R}\mathsf{Hom}(F(\mathcal{A}), F(\mathcal{B})) =
\mathsf{Hom}((F(\mathcal{A})_c,(F(\mathcal{B})_f)_{\lambda})\,.$$
where $\lambda$ denotes an infinite cardinal whose size is greater or
equal to the cardinality of the set of isomorphism classes of objects
in the category $\mathsf{H}^0(\mathcal{A}_c)$.
We dispose of the following $Q$-weak equivalent objects in $\mathsf{Lp}$:
\begin{eqnarray}
&&\mathcal{R}\mathsf{Hom}(F(\mathcal{A}), F(\mathcal{B})) \nonumber\\
&&\mathsf{Hom}((F(\mathcal{A})_c,(F(\mathcal{B})_f)_{\lambda})\nonumber\\
&&\mathsf{Hom}((\emptyset \subset \mathcal{A}_c), ((\mathcal{B}_f)_{\lambda})_{contr}
\subset (\mathcal{B}_f)_{\lambda}))\nonumber\\
&&(\mathsf{Fun}_{dg}(\mathcal{A}_c, ((\mathcal{B}_f)_{\lambda})_{contr}) \subset
\mathsf{Fun}_{dg}(\mathcal{A}_c, (\mathcal{B}_f)_{\lambda}))\nonumber\\
&&\overline{\mathsf{Fun}_{dg}(\mathcal{A}_c, ((\mathcal{B}_f)_{\lambda})_{contr})}
\subset \mathsf{Fun}_{dg}(\mathcal{A}_c,
(\mathcal{B}_f)_{\lambda}))/\mathsf{Fun}_{dg}(\mathcal{A}_c,
((\mathcal{B}_f)_{\lambda})_{contr})\nonumber\\
&&(\emptyset \subset \mathsf{rep}_{dg}(\mathcal{A}_c,\mathcal{B}_f))\nonumber\\
&&F(\mathsf{rep}_{dg}(\mathcal{A},\mathcal{B}))\nonumber\,.
\end{eqnarray}
This proves that the total right derived functor
$\mathcal{R}\mathsf{Hom}(-,-)$ in $\mathsf{Ho}(\mathsf{Lp})$
corresponds to the functor $\mathsf{rep}_{dg}(-,-)$, as in
\cite{dg-cat-survey} \cite{Toen}.  
\end{proof}

\end{document}